\documentclass[10pt]{article}

\usepackage[dvips]{graphicx}
\usepackage{amssymb}
\usepackage{amsmath}

\usepackage{caption}
\usepackage{hyperref}
\usepackage{arcs}
\usepackage {xcolor} 
\usepackage{subfigure}
\usepackage{authblk}

\newtheorem{theorem}{Theorem}

\newenvironment{AMS}{\small\bf 2010 AMS subject classification: }{} 

\begin{document}

\title{Numerical cubature on scattered data \\ by adaptive interpolation}

\author[1]{R. Cavoretto}
\affil[1,3]{University of Torino, Italy}
\author[2]{F. Dell'Accio}
\affil[2,4]{University of Calabria, Cosenza, Italy}
\author[3]{A. De Rossi}
\author[4]{F. Di Tommaso}
\author[5]{N. Siar}
\affil[5]{Faculty of Sciences, University Moulay Ismail-Meknes, Morocco}
\author[6]{A. Sommariva} 
\affil[6,7]{University of Padova, Italy}
\author[7]{M. Vianello}

\date{\today}
\maketitle

\begin{abstract} 
We construct cubature methods on scattered data 
via resampling on the support of known algebraic cubature formulas, by different kinds of adaptive interpolation (polynomial, RBF, PUM). This approach gives a promising alternative to other recent methods, such as direct meshless cubature by RBF or least-squares cubature formulas. 

\end{abstract}

\vskip0.2cm
\noindent
\begin{AMS}
{\rm 65D05, 65D32}
\end{AMS}

\section{Introduction}

The problem of computing integrals of multivariate functions on standard or nonstandard domains via function sampling is ubiquitous in scientific and technological applications. There is a vast literature on the construction of cubature formulas exact in specific function spaces, especially polynomial spaces, where the sampling has to be accomplished on suitable set of nodes. We do not even attempt to give an overview of this specialized literature, referring the reader to some classical or more recent monographs and surveys, such as for example \cite{C97,FG18,S71}.

On the other hand, it is a common situation that the data are scattered and still one would like to obtain an integration method from such data. One of the most popular approaches, especially in high-dimension, is QMC (Quasi MonteCarlo) method, provided that the sampling is made on low-discrepancy node sequences; again, we only mention within the vast QMC literature an excellent survey like \cite{DKS13}. 
While QMC has the advantage of not suffering from the curse of dimensionality, it is however scarcely accurate and thus not the best choice in low dimension.

Numerical cubature from scattered data that be exact on suitable function spaces, typically in low dimension, has received some attention in recent years. Indeed, we may quote methods based on meshless interpolation, namely by Radial Basis Functions (RBF), apparently considered for the first time in \cite{SV06,SV06-2} and then further developed in \cite{CDRSV22,GR23,SV21}. Another approach, based on $\ell^1$ or $\ell^2$ minimization of the weights under polynomial moment-matching conditions, was proposed in \cite{G21}, where it was shown to be more accurate than MC and QMC integration. 

Differently, here we construct cubature methods on scattered data via resampling on the support of known algebraic cubature formulas, by adaptive  interpolation of different kinds. The paper is organized as follows. In Section 2, we describe the fundamentals of our approach and the different kinds of interpolation adopted. In Section 3, we discuss a number of cubature tests, showing that the present method can be a valid alternative to the scattered integration methods quoted above.

\section{Cubature by adaptive scattered interpolation}
The basic idea of our approach is the following. Given an algebraic cubature rule with positive weights and interior nodes (PI rule) on a compact domain $\Omega\subset \mathbb{R}^s$
\begin{equation} \label{cformula}
I_\Omega(f)=\int_\Omega{f(P)\,w(P)\,dP}\approx Q_n(f)=\sum_{k=1}^\nu{w_k\,f(\Xi_k)}\;,
\end{equation}
which is exact for $f\in \mathbb{P}_n$ (the space of multivariate polynomials of total degree not exceeding $n$, with $m_s=dim(\mathbb{P}_n)={n+s\choose s}$), and a suitable interpolant on a scattered set $X=\{P_1,\dots,P_N\}\subset \Omega$, say $\Psi=\Psi_{f,X}:\Omega \rightarrow \mathbb{R}$, 
we approximate the cubature rule as
\begin{equation} \label{approx-cformula}
Q_n(f)=\sum_{k=1}^\nu{w_k\,f(\Xi_k)}\approx Q_n(\Psi)=\sum_{k=1}^\nu{w_k\,\Psi(\Xi_k)}\;.
\end{equation}
Notice that $\Psi$ could be a standard global interpolant such that $\Psi(P_i)=f(P_i)$, $1\leq i\leq N$, but also a function pointwise defined by local interpolation on a suitable subset of the scattered points, in a moving interpolation fashion    
(see for example Subsection 2.1). 

We stress that in this approach we do not approximate the integral of $f$ by the integral of the interpolant, which is the standard way to construct a cubature formula via interpolation, rather we consider the values of the interpolant at the cubature nodes of a known formula as perturbations of the function values there. In other words, we use the interpolant only as an approximate evaluator of the integrand at the cubature nodes, that is to {\em approximate the cubature formula}.

Now, denoting by $p^\ast_n\in \mathbb{P}_n$ the best uniform polynomial approximation to $f\in  C(\Omega)$, we can write
$$
I_\Omega(f)-Q_n(\Psi) = (I_\Omega(f)-I_\Omega(p^\ast_n))+(I_\Omega(p^\ast_n)-Q_n(p^\ast_n))
$$
\begin{equation} \label{split}
+(Q_n(p^\ast_n)-Q_n(f))+(Q_n(f)-Q_n(\Psi))
\end{equation}
and then, observing that the second summand is null by polynomial exactness of the cubature formula, 
we get 
$$
\left|I_{\Omega}(f)-Q_n(\Psi)\right| \leq 
I_\Omega(|f-p^\ast_n|)+Q_n(|p^\ast_n-f|)
+Q_n(|f-\Psi|)
$$
\begin{equation} \label{est2}
\leq \|w\|_{L^1} \left(2\,E_n(f;\Omega)+\max_k |f(\Xi_k)-\Psi(\Xi_k)|\right)\;,
\end{equation}
where $E_n(f;\Omega)=\inf_{p\in \mathbb{P}_n}\|f-p\|_{\infty,\Omega}$. 

We recall that the infinitesimal order of $E_n(f;\Omega)$ can be estimated on the so-called Jackson compact sets, that are compact domains for which there exist a positive integer $k_j$, $j=0,1,2,\dots$, and a constant $c_j(f)$ (depending on the partial derivatives of $f$ up to order $j$), such that if $f\in C^{k_j}(\Omega)$ a Jackson-like inequality holds, of the form $E_n(f)\leq c_j(f)\,n^{-j}$. Two basic examples are hypercubes with $k_j=j+1$  and Euclidean balls with $k_j=j$. In particular a sufficient
condition is that $\Omega$ be Whitney regular and admits a Markov polynomial inequality; cf. \cite{P09} for a survey of known results
on Jackson sets.

From (\ref{split})-(\ref{est2}) we expect that increasing $n$ the error of the approximate cubature $Q_n(\Psi)$ will stagnate on the size of the maximum interpolation error at the cubature nodes, since this term asymptotically dominates estimate (\ref{est2}). Such an error depends in turn on the choice of the interpolant and on the regularity of $f$. 

The reason for working with algebraic cubature is that Positive Interior (PI) formulas exact in polynomial spaces are available in a variety of domains. We cannot give here an extended bibliography, but we may quote some relevant items. For example, minimal or near-minimal formulas are known, at least up to certain degrees, on standard domains like squares/cubes, disks/balls and triangles/tetrahedra, cf. \cite{C97,CK00,FG18,FS12,S71,TFR22} with the references therein. On the other hand, low-cardinality formulas of PI type have been constructed on 
linear polygons/polyhedra \cite{JS20,SV23-2}, on domains with spline and NURBS boundaries \cite{SV23}, and on several domains corresponding to a trigonometric polynomial parametrization of the boundary such as for example circular segments, sectors, zones, lunes, as well as disk intersection and union, cf. \cite{DFSV13,DFV14,SV22} with the references therein. Many of these formulas use a NNLS implementation of Tchakaloff-like compression, cf. e.g. \cite{PSV17}. Other constructions of algebraic PI formulas on multidimensional domains can be found e.g. in \cite{KKN18,L21}, only to quote some recent contributions (see indeed the references therein). Thus, instead of trying to construct a cubature formula directly supported at the scattered sites, as done for example in \cite{CDRSV22,G21}, we use indirectly the scattered data to approximate one within such a vast and still growing body of available algebraic rules. 

As a preliminary study, below we briefly recall and study numerically several different adaptive interpolation methods on scattered data, comparing their approximation power on some test functions. Adaptivity is indeed the key to obtain high accuracy from such interpolants. We stress that here we are mainly interested in accuracy rather than computational efficiency, since 
our main goal will be approximating integrals as much accurately as possible by the fixed scattered data.

\subsection{Adaptive moving polynomial interpolation}
\label{sec.adamovpolint}

We discuss a pointwise evaluation method by adaptive local polynomial interpolation. 
In order to clarify its theoretical base, we state the main error estimate result, reformulating the theorem in \cite{DADTSV22}. Indeed, such a 
theorem concerns pointwise numerical differentiation of arbitrary order via local interpolation, whereas here we are merely interested in pointwise function evaluation. 
For simplicity, we restrict to bivariate functions. 

Below we denote by $\mathbb{P}_d$ the space of bivariate polynomials of total degree not exceeding $d$ with dimension 
$m_d=dim(\mathbb{P}_d)=(d+1)(d+2)/2$, 
by $B_h(\overline{P})$ the closed Euclidean ball
of radius $h$ centered at $\overline{P}=(\overline{x},\overline{y})$,
by $D^\alpha=\partial^{\alpha_1}_x\partial^{\alpha_2}_y$ the differentiation
operator with respect to bi-indexes 
$\alpha=(\alpha_1,\alpha_2)$ of length $|\alpha|=\alpha_1+
\alpha_2$ (lexicographically ordered), and by $C^{d,1}(K)$ the space
of $C^d$ functions with Lipschitz-continuous derivatives of
length $d$ on a compact domain $K\subset \mathbb{R}^2$ (the closure of a bounded open connected set), equipped with the
seminorm 
\begin{equation} \label{seminorm}
\|f\|_{C^{d,1}(K)}=\sup \left\{ \frac{|D^{\alpha }f(P)
-D^{\alpha }f(Q)|}{\|P-Q\|_2}:\,
P,Q\in K,\,P\neq Q,\,\left\vert \alpha \right\vert =d\right\}\;.
\end{equation}

Moreover, we adopt the usual notation with multi-indices where powers
are interpreted componentwise, $P^{\alpha}=x^{\alpha_1}y^{\alpha_2}$.

\begin{theorem}
Let $\Omega\subset \mathbb{R}^2$ be a convex body,
$f\in C^{d,1}(\Omega)$, and
$\pi_{h,d}
\in \mathbb{P}_d$ the interpolating polynomial of $f$ at a unisolvent
subset
\begin{equation} \label{Xm}
X_d=\{\xi_1,\dots,\xi_{m_d}\}
\subset \mathcal{N}_h=B_h(\overline{P})\cap \Omega\;.
\end{equation} 
Then the following pointwise error estimate holds
\begin{equation} \label{pointwise}
E_{h,d}(\overline{P})=\left|f(\overline{P})-\pi_{h,d}(\overline{P})\right|\leq \lambda_{h,d}(\overline{P})\,C_{h,d}\,h^{d+1}\;,
\end{equation}
$$C_{h,d}=\frac{2^d}{(d-1)!}\,\|f\|_{C^{d,1}(\mathcal{N}_h)}\;,\;\; \lambda_{h,d}(\overline{P})=\|\rho_1(V_{h,d}^{-1})\|_1\;,$$
where $\rho_1$ denotes the first row of the inverse Vandermonde matrix in the scaled monomial basis centered at $\overline{P}$
\begin{equation} \label{vand} 
V_{h,d}=\left[((\xi_i-\overline{P})/h)^\alpha\right]\;,\;1\leq i\leq m_d\;,\;|\alpha|\leq d\;.
\end{equation}

\end{theorem} 
\vskip0.5cm 

For a proof of this result, based on local approximation by Taylor formula, we refer the reader to \cite{DADTSV22}. 
We stress that $\lambda_{h,d}(\overline{P})$ is nothing but the Lebesgue function of $X_d$ evaluated at the neighborhood of the center $\overline{P}$, where the representation as 1-norm of the first Vandermonde row comes from the choice of the polynomial basis.
The result above can be extended to non-convex domains that satisfy the well-known ``Whitney regularity'' property \cite[\S 2]{W34}, by a ``curved'' version of Taylor formula, cf. e.g. \cite[Ch.8]{MP20}. In this case estimate (\ref{pointwise}) is qualitatively the same, since the Taylor formula remainder is still $\mathcal{O}(h^{d+1})$ (we do not give a proof for the sake of concision).  

The algorithm here adopted (2-dimensional moving interpolation), which is a derivation of the algorithm DISC proposed in \cite{DADTSV22-2} where we refer the reader for all details, works adaptively on the radius $h$ and the degree $d$, trying to minimize (an estimate of) the pointwise evaluation error. The local interpolation points $X_d$ are ``discrete Leja points'' extracted from $X\cap \mathcal{N}_h$ (we recall that $X$ are the fixed scattered sites), that are points aimed at maximizing the Vandermonde determinant modulus. 
They are selected greedily by standard Gaussian elimination with row pivoting on the local transposed Vandermonde matrix, following the method of \cite{BDMSV10} for stable multivariate polynomial interpolation.

The choice of discrete Leja points as local interpolation points is suggested by several features. First, being aimed at maximizing (in a greedy way) the Vandermonde determinant modulus, they guarantee {\em unisolvency}, provided that the local scattered points $X\cap \mathcal{N}_h$ are $\mathbb{P}_d$-determining (i.e., polynomials in $\mathbb{P}_d$ vanishing there vanish everywhere, or equivalently the corresponding rectangular Vandermonde matrix is full-rank). Though this is not true in general with any scattered sample, it is almost surely verified with uniformly distributed points with respect to any continuous density, as recently proved in \cite{DASV23} for general interpolation by a.e. analytic functions.

The second reason for the choice of discrete Leja points is interpolation {\em stability}. Indeed, their choice tries to keep small the entries of the inverse Vandermonde matrix (which are cofactors divided by the
Vandermonde determinant), and thus also the relevant row 1-norm in estimate (\ref{pointwise}). We observe also that the choice of a monomial basis centered at $\overline{P}$, and scaled by the ball radius, allows to control the Vandermonde matrix conditioning, at least up to moderate interpolation degrees (say around degree 10).

The third important feature of discrete Leja points is that they form a {\em sequence}, differently from other extremal interpolation sets like discrete Fekete points \cite{BDMSV10}. This means that the first $m_k$ among $m_d$ Leja points are unisolvent for degree $k=1,2,\dots,d$, a feature that is exploited by the algorithm DISC in \cite{DADTSV22-2} to compute a pointwise a posteriori error estimate.

\begin{figure}
\centering
\includegraphics[height=1.7in]{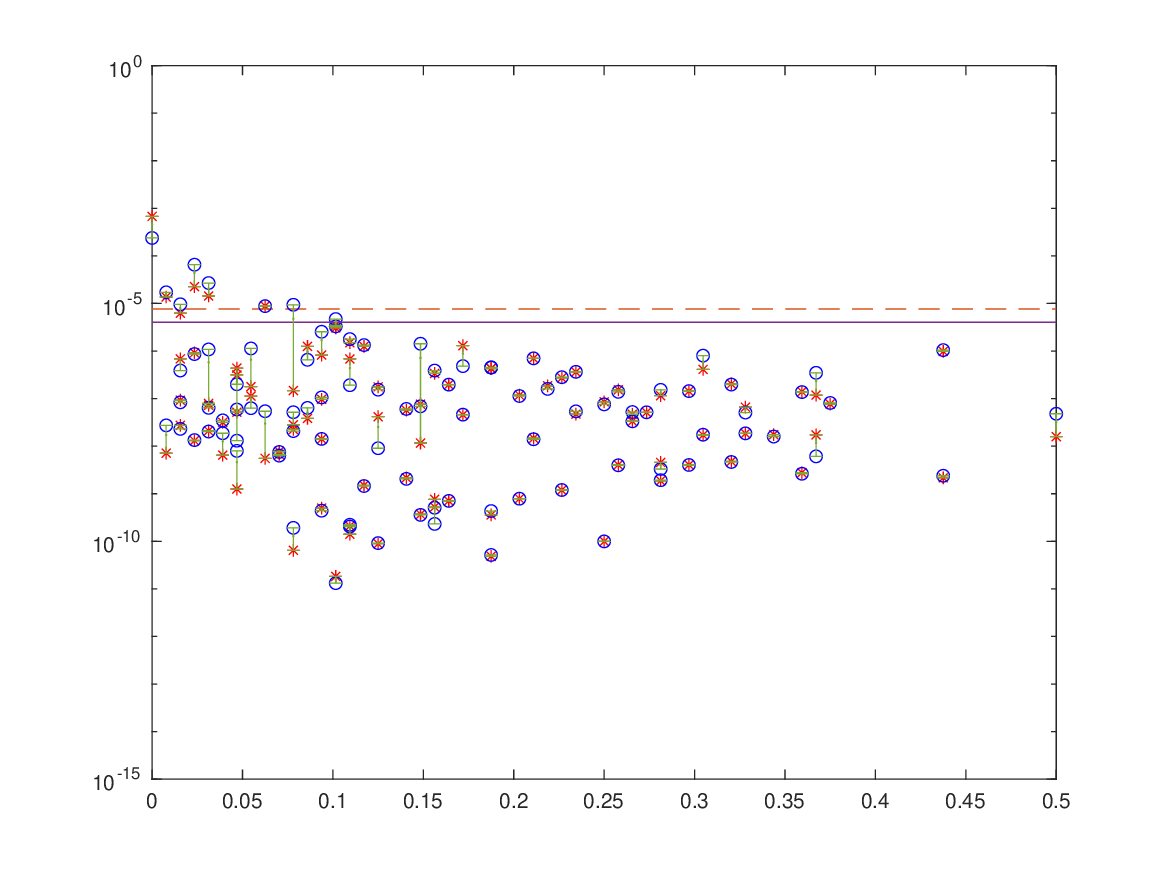}
\includegraphics[height=1.7in]{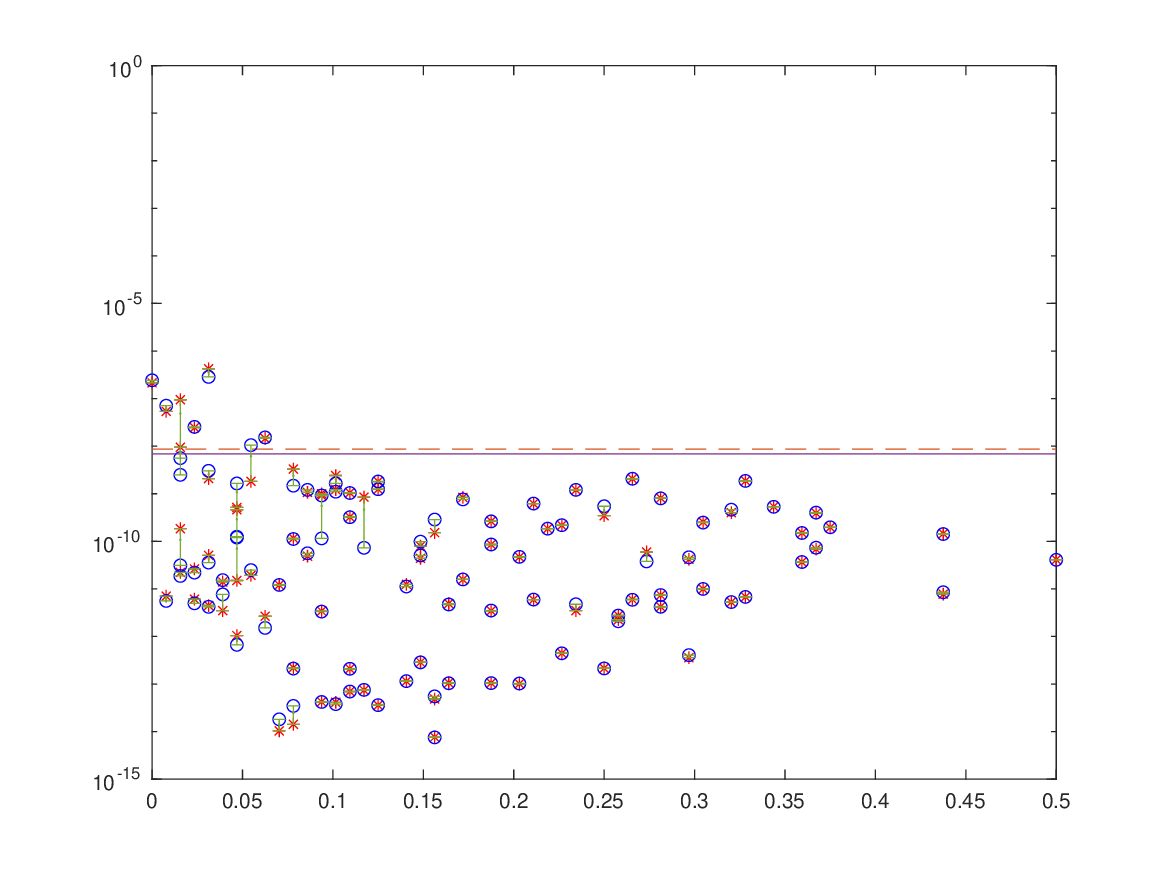}
\caption{Pointwise evaluation errors  (circles) and estimates (asterisks) by adaptive moving interpolation of Franke’s test function on 100 Sobol test 
points in $[0,1]^2$, ordered
by distance from the boundary, using 800 (left) and 1600 (right) 
Halton interpolation points (solid line: mean error; dashed line: mean estimate).}
\label{Franke}
\end{figure}

In order to show the performance of moving interpolation, in Figure \ref{Franke} we report the pointwise error and the error estimate in \cite{DADTSV22-2}  for the well-known Franke's test function, by moving interpolation on 800 and 1600 Halton points in $[0,1]^2$. These quantities are computed for 100 Sobol evaluation points, ordered by increasing distance from the boundary. Observe that in most cases the estimate is close to the actual error and only in few cases it substantially overerestimates or
underestimates the error, by one or seldom two orders of magnitude. The resulting mean estimate (dashed line) turns out to be a good approximation of the mean error (solid line); as expected, both decrease by increasing the number of sampling points.
Notice that the errors above the mean tend to occur near the boundary. This boundary effect, already observed in \cite{DADTSV22,DADTSV22-2}, can be explained observing that near the boundary the neighborhoods $\mathcal{N}_h=B_h(\overline{P})\cap \Omega$ contain less sampling points.  
Consequently, the maximum local interpolation degree becomes smaller and at the same time the quality of the Leja-like interpolation points can be worse.

\subsection{Adaptive LOOCV RBF interpolation} \label{adapLOOCVRBF}

We now provide a brief overview on adaptive LOOCV RBF interpolation recalling basic notions on RBF theory, which is useful for scattered data approximation; for further details, see \cite{fas15,wen05}.	

We suppose that is given a compact domain $\Omega \subset \mathbb{R}^s$, a scattered point set $X=\{P_1,\dots,P_N\}\subset \Omega$, and the corresponding data (or function) value set $F=\{f_1,\ldots,f_N\}\subset \mathbb{R}$ that is obtained by possibly sampling any (unknown) function $f: \Omega \rightarrow \mathbb{R}$. 

Considering a RBF $\phi: [0,\infty)\to \mathbb{R}$ that is \emph{strictly conditionally positive definite} of order $m$ (SCPD$_m$) \cite{fas07}, if we set $\phi_i(P)=\phi(||P-P_i||_2)$, one can find a unique interpolating function $\Psi: \Omega \rightarrow \mathbb{R}$ of the form
\begin{equation} \label{rbfi}
	\Psi(P) = \sum_{i=1}^{N} c_i \phi_i (P) + \sum_{i=N+1}^{N+M} c_{i}\pi_{i-N}(P),
\end{equation} 
where $\{\pi_k\}_{k=1}^M$ generate a basis for the 	$M = {m-1+s \choose m-1}$-dimensional linear space $\mathbb{P}_{m-1}^s$ of $s$-variate real valued polynomials of total degree less than or equal to $m-1$, and $||\cdot||_2$ is the Euclidean norm. The coefficients $c_1,\ldots,c_{N+M}$ are determined by enforcing the interpolation conditions $\Psi(P_i)=f_i$, for $i=1, \ldots, N$. Since these conditions lead to a system of $N$ linear equations in the $N+M$ unknowns $c_i$, one usually adds the $M$ additional conditions $\sum_{i=1}^{N} c_{i}\pi_{k}(P_{i})=0$, $k=1,\dots, M$, to ensure a unique solution. Moreover, from theory it is known that a SCPD$_0$ function is \emph{strictly positive definite} (SPD), and so in this case the polynomial in (\ref{rbfi}) is omitted \cite{fas07} 

Solving the interpolation problem for a SCPD$_m$ function $\phi$ leads to a symmetric linear system
\begin{equation} \label{mat}
{\cal M}\textbf{c}=\textbf{y}, 
\end{equation}
where
$$
{\cal M}=\left[
\begin{array}{cc}
A   & B  \\
B^T & O 
\end{array}
\right], \qquad
\textbf{y}=\left[
\begin{array}{c}
\textbf{f}  \\
\textbf{0} 
\end{array}
\right].
$$
The \emph{interpolation matrix} ${\cal M}$ in (\ref{mat}) has entries $A_{ki}=\phi(||P_k-P_i||_2)$, $B_{kj}=\pi_j(P_k)$, with $k,i=1,\ldots,N$, $j=1,\ldots,M$, and $O$ is a $M \times M$ zero matrix. Moreover, $\textbf{c}=[c_1,\ldots,c_{N+M}]^T$, $\textbf{f}=[f_1,\ldots,f_N]^T$ and $\textbf{0}$ is a zero vector of length $M$. Notice that for a SPD function $\phi$ the matrix reduces simply to ${\cal M} = A$, and the polynomial part vanishes.

In literature, many RBFs are scaled by a \textsl{shape parameter} $\varepsilon > 0$ such that 
$$
\phi_{\varepsilon}(r)=\phi(\varepsilon r)=\phi(\varepsilon ||P-P_i||_2), \quad \forall P, P_i \in \Omega. 
$$
Hereinafter, for the sake of simplicity, we keep implicit (unless strictly necessary) the dependence on $\varepsilon$, referring to $\phi_i(P)$ as $\phi_{\varepsilon}(r)$. For such RBFs the choice of a \lq\lq good\rq\rq, or possibly an \lq\lq optimal\rq\rq, shape parameter is a crucial task, but also a big issue (see e.g. \cite{cav21a} and references therein). Some examples of popular SCPD RBFs together with their smoothness and abbreviation are listed as follows \cite{fas07}: 
$$
\phi_{\varepsilon}(r) =\left\{
\begin{array}{lll}
\exp(-\varepsilon^2 r^2), & & \quad \mbox{Gaussian $C^{\infty}$} (\mbox{GA})  \medskip  \\ 
(1+\varepsilon^2r^2)^{-1/2}, & & \quad \mbox{Inverse MultiQuadric $C^{\infty}$} (\mbox{IMQ})  \medskip  \\
(1+\varepsilon^2r^2)^{1/2}, & & \quad \mbox{MultiQuadric $C^{\infty}$} (\mbox{MQ})  \medskip  \\
\max \left(1-\varepsilon r,0\right)^4(4\varepsilon r+1), & & \quad \mbox{Wendland $C^2$} (\mbox{W2}) 
\end{array}
\right.
$$
All these functions are SPD, except for MQ that is SCPD$_1$. GA, IMQ and MQ are globally supported, while W2 is compactly supported and its support is $\left[0,1/\varepsilon\right]$, see \cite{wen05}.

However, as remarked in \cite[Ch.11]{fas15}, from a correct reformulation of the so-called \emph{uncertainty} or \emph{trade-off} \emph{principle} due to Schaback \cite{sch95a} we know that using a standard basis one cannot have high accuracy and stability at the same time. Indeed, when interpolating target functions by RBFs, the best accuracy is typically achieved in the \emph{flat limit} $\varepsilon \rightarrow 0$, but in this case the interpolation matrix might be severely ill-conditioning. Thus, in order to get trustworthy results, in the literature several techniques have been proposed to guide us in a suitable selection of the RBF shape parameter (see e.g. \cite[Ch.14]{fas15} and \cite{gol15,sch11}).

A popular strategy for estimating the RBF shape parameter $\varepsilon$ is the LOOCV method. It was originally introduced in RBF interpolation by Rippa \cite{rip99}, and more lately has widely been used and extended in many other fields (see e.g. \cite{cav20a,cav20b,lin22}). 

In the LOOCV technique the shape parameter is selected by minimizing a cost function that collects the errors for a sequence of partial fits to the data. To estimate the unknown true error, we split the data into two parts: a \emph{training} data set consisting of $N-1$ data to obtain a \lq\lq partial interpolation\rq\rq, and a \emph{validation} data set that contains a single (remaining) data used to compute the error. After repeating in turn this procedure for each of the $N$ given data, the result is a vector of error estimates and the cost function is used to optimally determine the value of $\varepsilon$.

For this discussion we define by $X^{[k]}=\{P_1,\ldots,P_{k-1},P_{k+1},\ldots,P_N\}$, and $F^{[k]}=\{f_1,\ldots,f_{k-1},f_{k+1},\ldots,f_N\}$, the sets of points and corresponding values with the removed data $(P_{k},f_k)$, denoted by the superscript $[k]$. All other quantities are represented similarly.  

The key idea of LOOCV is to predict the parameter $\varepsilon$ by the partial RBF interpolant to the data $(X^{[k]},F^{[k]})$, i.e.,
\begin{equation} \label{pint}
	\Psi^{[k]}(P) = \sum_{i=1,\ i \neq k}^{N} c_i \phi_i(P) + \sum_{i=N+1}^{N+M} c_i\pi_{i-N}(P).
\end{equation}
The interpolant (\ref{pint}) enables us to find the error $e_k(\varepsilon) = f_k - \Psi^{[k]}(P_k)$, $k=1,\ldots,N$, removing in turn each point $P_k$ and comparing then the resulting fit with the known value at the removed point $f_k$. This error can be computed more efficiently by simply using the rule
\begin{equation} \label{ekeps} 
e_k(\varepsilon) = \frac{c_k}{{\cal M}^{-1}_{kk}},
\end{equation}
where $c_k$ is the $k$th coefficient deriving from the solution of the \emph{full} RBF system (\ref{mat}), and ${\cal M}^{-1}_{kk}$ is the $k$th diagonal element of the matrix ${\cal M}^{-1}$. This formulation needs to only solve a single linear system, considering the entire data set $(X,F)$.  

In summary, the problem is solved by minimizing the LOOCV cost function
$$
\mbox{LOOCV}(\varepsilon) = ||\textbf{e}(\varepsilon)||_{\infty} = \max_{k=1,\ldots,N} \left|\frac{c_k}{{\cal M}^{-1}_{kk}}\right|.
$$
where $\textbf{e}(\varepsilon)=\left[e_k(\varepsilon)\right]_{k=1,\dots,N}$ and $||\cdot||_{\infty}$ denotes the $\infty$-norm, even if in principle any norm is enabled.



\subsection{BLOOCV RBF-PUM}

We now introduce the basic theory for an adaptive RBF-PUM interpolation, which is based on the BLOOCV scheme \cite{CDRP18}. Since inverting the kernel matrix in (\ref{mat}) is computationally expensive when the amount of data highly increases, an effective way to overcome this issue is to partition the open bounded domain $\Omega$ into $p$ overlapping subdomains $\{\Omega_{j}\}_{j=1}^{p}$ which form a open bounded cover of $\Omega$, that is $\Omega \subseteq \bigcup_{j=1}^p\Omega_j$. So the interpolation problem is locally decomposed into $p$ (smaller) subproblems, where the subdomains of the PU covering consist of overlapping balls of radius $\delta_j$, see \cite{cav21}. 

The idea of PUM comes originally from the context of PDEs \cite{Babuska97}, then gaining popularity in the RBF community in various fields of applied mathematics (see e.g. \cite{cav12,lar17,wen02}). The PUM has indeed some features that makes it particularly suitable for the approximation of large and irregularly distributed scattered data sets. In particular, it enables us to efficiently solve a local interpolation problem on each PU subdomain and thus to construct the global approximant by gluing together the local contributions using some weight functions. To achieve that, we need those weights to be a family of compactly supported, non-negative, and continuous functions $w_j:\Omega_j \rightarrow \mathbb{R}$, with $\mbox{supp}(w_j) \subseteq \Omega_j$, such that
$$
\sum_{j=1}^{p} w_j\left(P\right) = 1, \quad P \in \Omega.
$$
Once we choose the partition of unity $\{ w_j \}_{j=1}^{p}$, we form the global interpolant 
\begin{equation} \label{intg}
\Psi(P) = \sum_{j=1}^{p} \Psi_j(P) w_j(P), \quad P \in \Omega,
\end{equation}
which is a weighted sum of $p$ local RBF interpolants $\Psi_j:\Omega_j \rightarrow \mathbb{R}$, i.e.
\begin{equation} \label{intl}
\Psi_j(P) = \sum_{i=1}^{N_j} c_i^j \phi_i^j(P),
\end{equation}
where $\phi_i^j(P) = \phi_{\varepsilon_j}(||P - P^j_i||_2)$, $P_i^j \in X_j = X \cap \Omega_j$, and $N_j=|\Omega_j|$. The interpolant (\ref{intl}) represents the local version of (\ref{rbfi}) in the special case of SPD RBFs, and it is uniquely found by solving the linear system associated with the subdomain $\Omega_j$, i.e.,
$$
A_j\textbf{c}_j = \textbf{f}_j,
$$
where $(A_j)_{ki}=\phi_{\varepsilon_j}(||P_k^j - P_i^j||_2)$, $k,i=1,\ldots,N_j$, $\textbf{c}_j=[c_1^j,\ldots,c_{N_j}^j]^T$, $\textbf{f}_j=[f_1^j,\ldots,f_{N_j}^j]^T$; for a generalization to SCPD$_m$ functions, see \cite{wen02}. Moreover, from \cite{CDRP18} it is known that the accuracy of the global fit (\ref{intg}) strongly depends on the choices of the shape parameter $\varepsilon_j$ and the radius $\delta_j$.

The LOOCV technique discussed in Subsection \ref{adapLOOCVRBF} can thus be extended to search the optimal values of $\varepsilon_j$ and $\delta_j$ in each PU subdomain. It consists in evaluating in every subdomain $\Omega_j$, for each couple of $(\varepsilon,\delta)$ and $k \in \lbrace 1, \dots, N_j \rbrace$, the interpolation error $e^j_k(\varepsilon,\delta)$ at the point $P_k$ of the partial RBF interpolant $\Psi_j^{[k]}$ fitted on the point set $X_j^{[k]} = X_j \setminus \lbrace P_k^j \rbrace$ and the data value set $F_j^{[k]} = F_j \setminus \lbrace f_k^j \rbrace$. It consists in computing in every subdomain $\Omega_j$, for each couple of $(\varepsilon,\delta)$ and $k \in \lbrace 1, \dots, N_j \rbrace$, the interpolation error $e^j_k(\varepsilon,\delta)$  of the partial RBF interpolant $\Psi_j^{[k]}$ evaluated at the point $P^j_k$, fitted on the point set $X_j^{[k]} = X_j \setminus \lbrace P_k^j \rbrace$ and the data value set $F_j^{[k]} = F_j \setminus \lbrace f_k^j \rbrace$. To avoid this computation, by suitably adapting (\ref{ekeps}), the local error is defined as follows:
$$
e_k^j(\varepsilon,\delta) = \frac{c_k^j}{(A_j^{-1})_{kk}}.
$$

Hence, the optimal value of $(\varepsilon_j,\delta_j)$ is the one that minimizes the error function 
$$
\mbox{LOOCV}_j(\varepsilon,\delta) = \max_{k=1, \dots, N_j} \left|   \frac{c_k^j}{(A_j^{-1})_{kk}} \right|.
$$

\subsection{Multinode Shepard method} 
 The Multinode Shepard (MS) method is an accurate procedure for reconstructing functions from scattered data. It is the generalization of the Little's idea that consists in improving the classic Shepard method by acting in two ways:  modifying
the classical point-based weight functions and defining instead a normalized blend of locally
linear interpolating polynomials with triangle-based weight functions which depend on the product
of inverse distances to the three vertices of the corresponding triangle \cite{DADTH15}. We assume that $\Omega \subset \mathbb{R}^s$, $s\geq 2$, is a region, $X=\{P_i\}_{i=1}^{N}\subset\Omega$ is a finite set of pairwise distinct scattered points and $F=\{f_i\}_{i=1}^N$ is a set of associated real numbers, obtained by sampling a continuous function $f:\Omega\rightarrow {\mathbb R}$. By working with local interpolating polynomials of total degree $d\in \mathbb{N}$ of $s$ variables, we let $m_d=\binom{s+d}{d}$ the dimension of such a polynomial space. To define the multinode Shepard method, we need the existence of a set $\mathcal{S}=\{\sigma_j\}_{j=1}^{S} \subset X$, which is a covering of the set $X$, i.e. 
$$
 \bigcup_{j=1}^S \sigma_{j}= X,
$$
 and such that each $\sigma_j=\{P_{j_1},\dots,P_{j_{m_d}}\}$ is unisolvent for polynomial interpolation of total degree $d$. As proven in \cite{DASV23}, such a set $\mathcal{S}$ exists for almost all choices of the interpolation nodes.  The local interpolation polynomial $\pi_j(P)$ based on $\sigma_j=\{P_{j_1},\dots,P_{j_{m_d}}\}$ is represented as 
\begin{equation*}
    \pi_j(P)=\sum_{i=1}^{m_d} \ell _{j,i}(P) f_{j_i}\;,\;\;P \in \mathbb{R}^s,
\end{equation*}
 where
 \begin{equation*}
  \ell _{j,i}\left( P\right)=   \sum\limits_{\underset{|\alpha|\leq d}{\alpha \in \mathbb{N}^s_{0}}} a_{\alpha}^{(i)}\left(P-P_{j}^{(b)}\right)^{\alpha},\;\;P \in \mathbb{R}^s,
\label{lagrange_basis_functions}
 \end{equation*}
 are the Lagrange basis polynomials satisfying
 $$
 \ell _{j,i}(P_{j_k})=\delta_{ik}=
 \left\{
 \begin{array}{rl}
 1, & i=k,\\
 0, & i\ne k.
 \end{array}
 \right.
 $$
 written, using a multi-index notation, in the Taylor basis centered at the barycenter $P_{j}^{(b)}$ of $\sigma_j$ \cite{DADTS21}. Likewise the triangular \cite{DADTH15} ($s=2$, $d=1$),  the hexagonal \cite{DADT20} ($s=2$, $d=2$) and the tetrahedral Shepard methods \cite{CDRDADT20} ($s=3$, $d=1$), the MS method is obtained by combining  the polynomials $\ell_{j,i}\left( P\right)$, $i=1,\dots,N,\ j=1,\dots,S$ with the MS basis function \cite{DADT19}
\begin{equation}
W_{\mu ,j}\left( P\right) =\dfrac{\prod\limits_{l=1}^{m_d}\left\Vert 
P-P_{j_{l}}\right\Vert_{2} ^{-\mu }}{\sum\limits_{k=1}^{S}%
\prod\limits_{l=1}^{m_d}\left\Vert P-P_{k_{l}}\right\Vert_{2}
^{-\mu }}\;,\;\;\mu>0\;, \hspace{0.1cm}P \in \mathbb{R}^s,
\label{m_points_basis_functions}
\end{equation}%
as follows
\begin{equation}
\Psi_{\mu }\left( P\right) 
=\sum\limits_{i=1}^{N}\sum\limits_{j\in \mathcal{J}_{i}}W_{\mu ,j}\left( 
P\right) \ell _{j,i}\left( P\right) f_i\;,
\label{multinode_rewriting}
\end{equation}
where $\mathcal{J}_{i}=\{j\in \{1,\dots ,S\}: P_i\in \sigma_j \}$. The functions $W_{\mu,j}\left(P\right)$ are a partition of unity and the MS approximant interpolates on all the scattered points $P_i$, $i=1,\dots,N$ and reproduces polynomial of degree less than or equal to $d$.

Results on the approximation order of the MS method can be found in \cite{DADT19}. By letting $||\cdot ||_\infty$
the maximum norm, $R_{r}(Q)=\{P\in \mathbb{%
\mathbb{R}
}^{s}:||P-Q||_\infty\leq r\}$ the axis-aligned closed cube with center $Q$ and edge
length $2r$, we set%
\begin{equation}
h^{\prime }=\inf \{r>0:\forall P\in \Omega \ \exists \sigma_j\in \mathcal{S}:R_{r}(P)\cap
\sigma_j\neq \emptyset \},  \label{eq:h-prime}
\end{equation}%
\begin{equation}
h^{\prime \prime }=\inf \{r>0:\forall \sigma_j\in \mathcal{S}\ \exists P\in \Omega
:\sigma_j\subset R_{r}(P)\}  \label{eq:h-second}
\end{equation}%
and finally
\begin{equation}
h=\max \{h^{\prime },h^{\prime \prime }\}.  \label{eq:h}
\end{equation}%
The positive real number $h$ is clearly a measure of the fill distance of the 
points in $X$ and of the largeness of subsets' $\sigma_j$ diameters: $h$ decreases if
the number of a rather uniform distribution of scattered points increases
and the diameters of the subsets $\sigma_j$ remain relatively small. We further let
\begin{equation} 
M=\sup_{P\in \Omega }\#\{\sigma_j\in \mathcal{S}:R_{h}(P)\cap \sigma_j\neq \emptyset \},
\label{eq:magicM}
\end{equation}%
the maximum number of subsets $\sigma_j$ with at least one point in some square with
edge length $2h$. Small values of $M$, in correspondence of small values of $%
h$, imply that there are no clusters of subsets of the covering $\mathcal{S}$. Under the hypothesis of compactness and convexity of the region $\Omega $, if $\mu >\frac{s+d+1}{m_d}$ in \cite{DADT19} it is proven that 
\[
||f-\Psi_{\mu }\left( P\right)||_{\infty}\leq CM||f||_{d,1}h^{d+1},
\]%
where $C$ is a positive constant and 
$$
||f||_{d,1}=\sup\left\{\dfrac{\left|D^{\nu}f(P)-D^{\nu}f(Q)\right|}{||P-Q||_{2}},\,P,Q \in \Omega,\,\nu\in\mathbb{N}^s_{0},\,\left|\nu\right|=d\right\},
$$
provided that $f$ is differentiable with partial derivatives of order $d$ Lipschitz-continuous. The constant $C$ is explicitly computed; it is directly proportional to $\left(1+\max\limits_{j=1,\dots,S}||\pi_j||_\infty\right)$ by a factor which depends only on $s$, $S$, $d$, $m_d$ and $\mu$. To maintain $\max\limits_{j=1,\dots,S}||\pi_j||_\infty$ small, \cite{DADTNS21}, a convenient procedure consists in considering, for each scattered point, the set of $m_d+q$, $q>0$, nearest points and in choosing, among them, the subset of $m_d$ discrete Leja points by the procedure introduced in \cite{BDMSV10} detailed in section \ref{sec.adamovpolint}. Finally, the result on the order of convergence holds with a slightly different constant in the case of non-convex domains (even multiply connected domains)  by using the ``Whitney regularity'' property, in line with what is observed in section \ref{sec.adamovpolint} (see also \cite{F986}).

\begin{figure}
\centering
\includegraphics[scale=0.37]{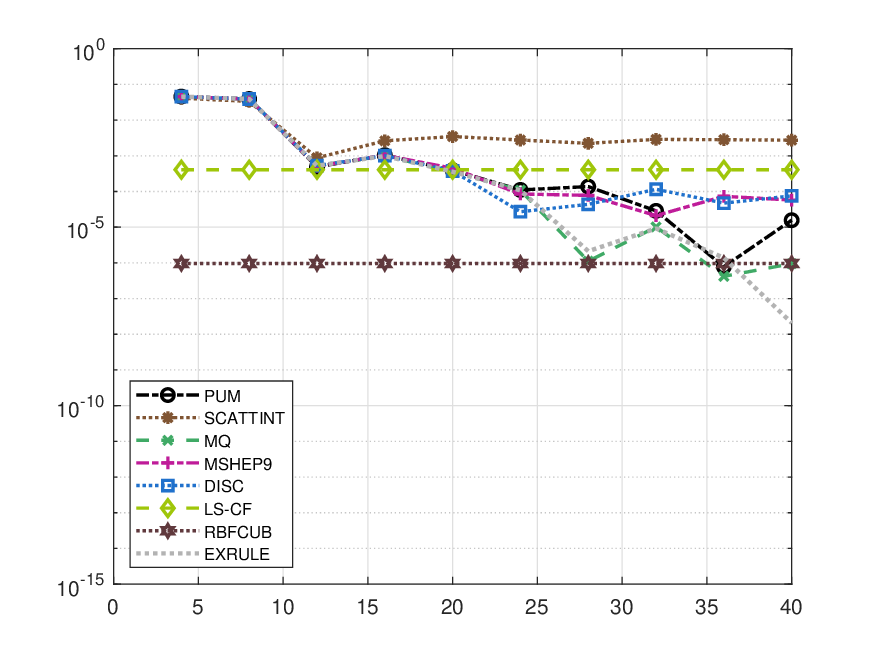}
\includegraphics[scale=0.37]{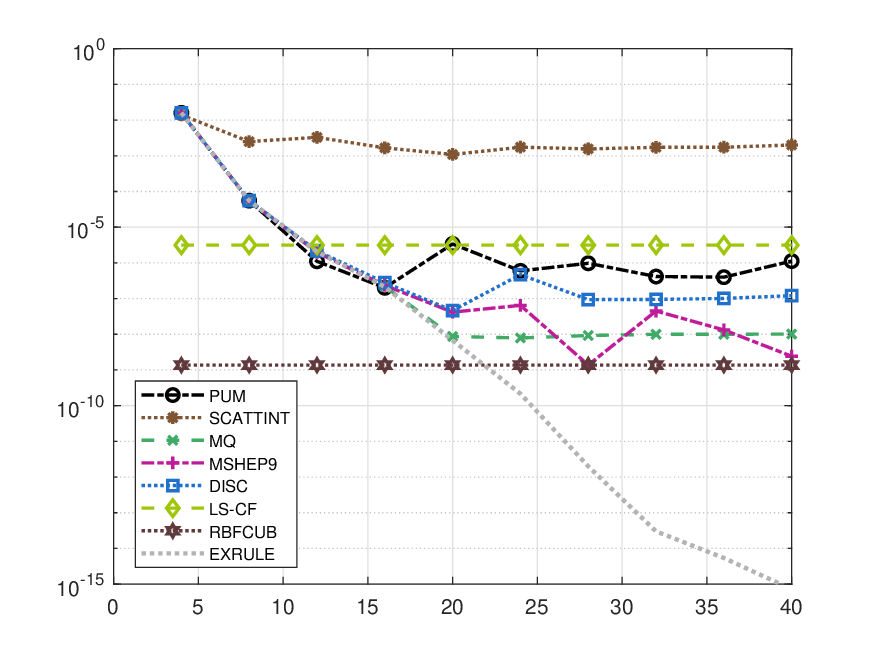}
\includegraphics[scale=0.37]{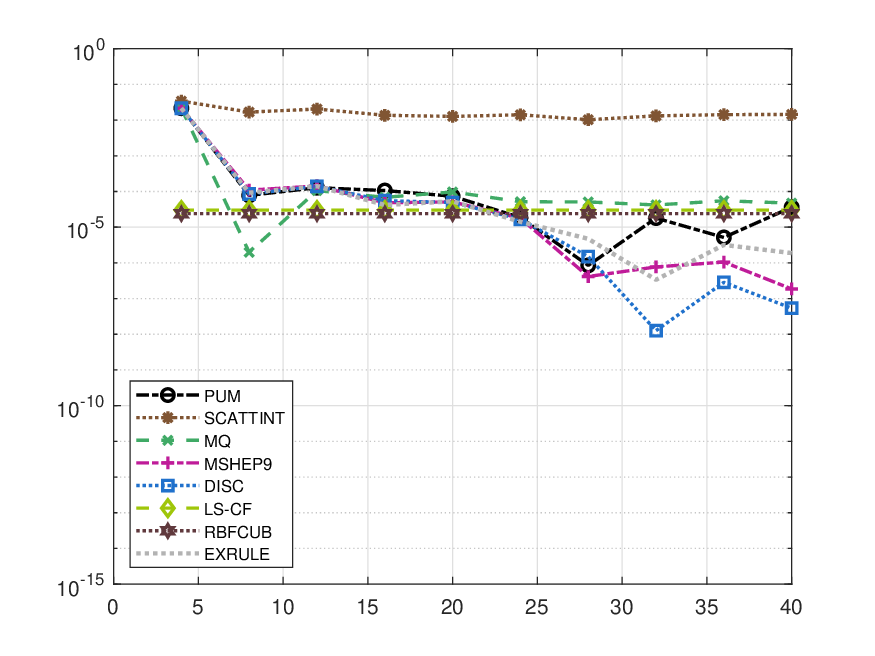}
\includegraphics[scale=0.37]{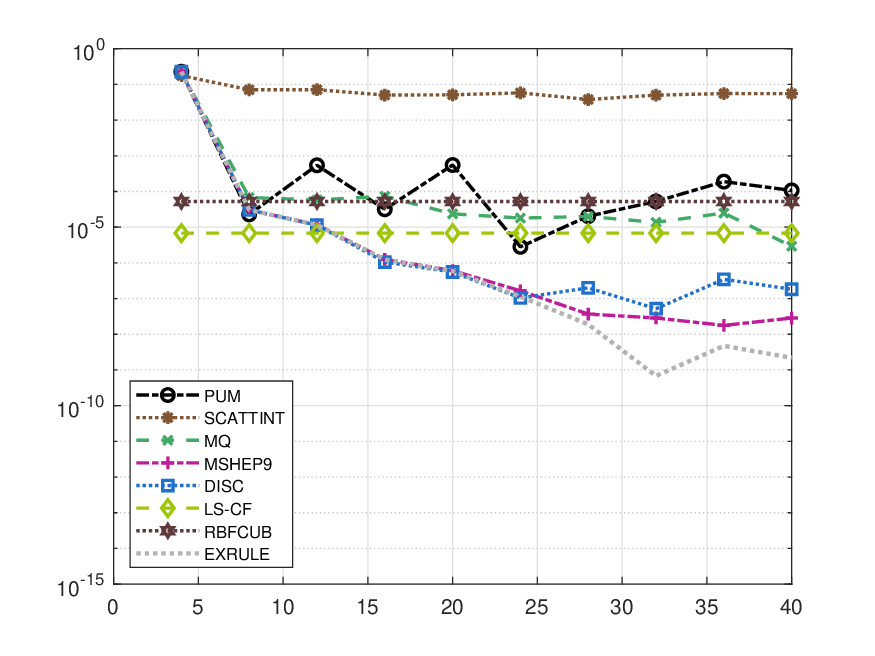}  
\caption{Comparison of scattered cubature methods on 400 Halton points for the test functions $f_1$ (top left) and $f_2$ (top right), $f_3$ (bottom left) and $f_4$ (bottom right); the horizontal lines correspond to RBFCUB with MultiQuadrics and to LS-CF method, whereas the grey dotted line to the underlying algebraic rule with exact function values.}
\label{test0}
\end{figure}

\section{Cubature tests}

In this section we compare cubature by adaptive scattered interpolation via suitable PI algebraic rules as in (\ref{cformula})-(\ref{approx-cformula}), with two relevant methods proposed in the recent literature, namely RBFCUB in \cite{CDRSV22} (based directly on LOOCV optimized RBF integration over polygons) and Glaubitz's LS-CF in \cite{G21} (based on $\ell^2$ weight minimization under polynomial moment matching conditions). To our knowledge, this is the first systematic comparison of different approaches for numerical cubature by scattered data, in particular concerning how much accuracy is obtainable from the data. All the numerical tests have been performed by the Matlab package SCATTCUB, freely available along with the corresponding demos at \cite{CDRDADTSSV23}.

\begin{figure}
\centering
\includegraphics[scale=0.37]{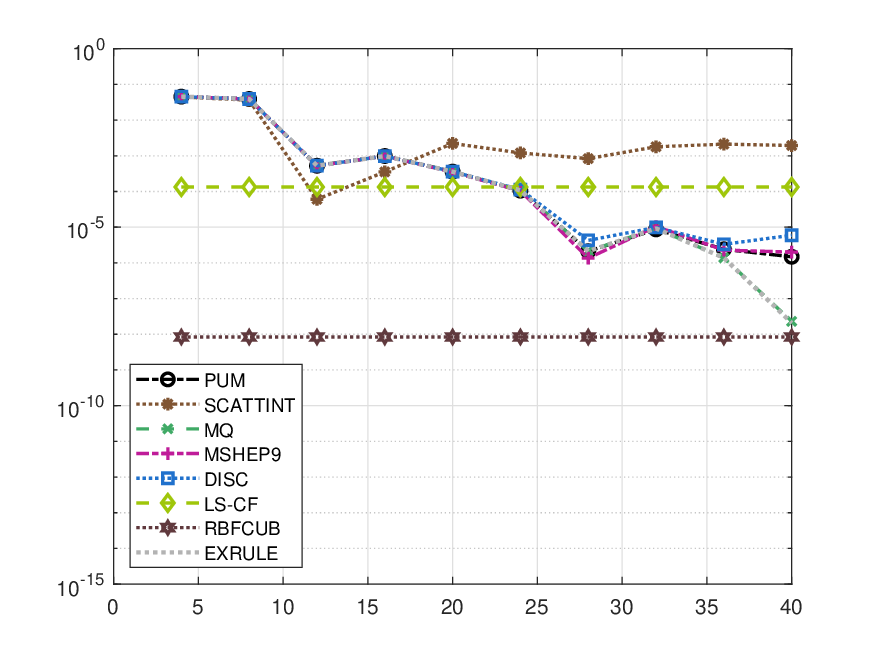}
\includegraphics[scale=0.37]{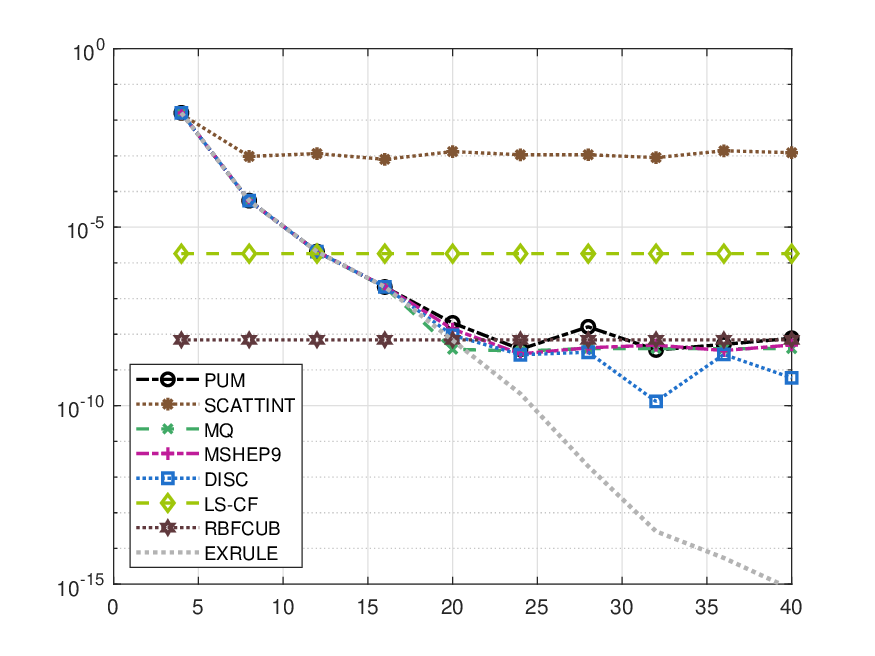}
\includegraphics[scale=0.37]{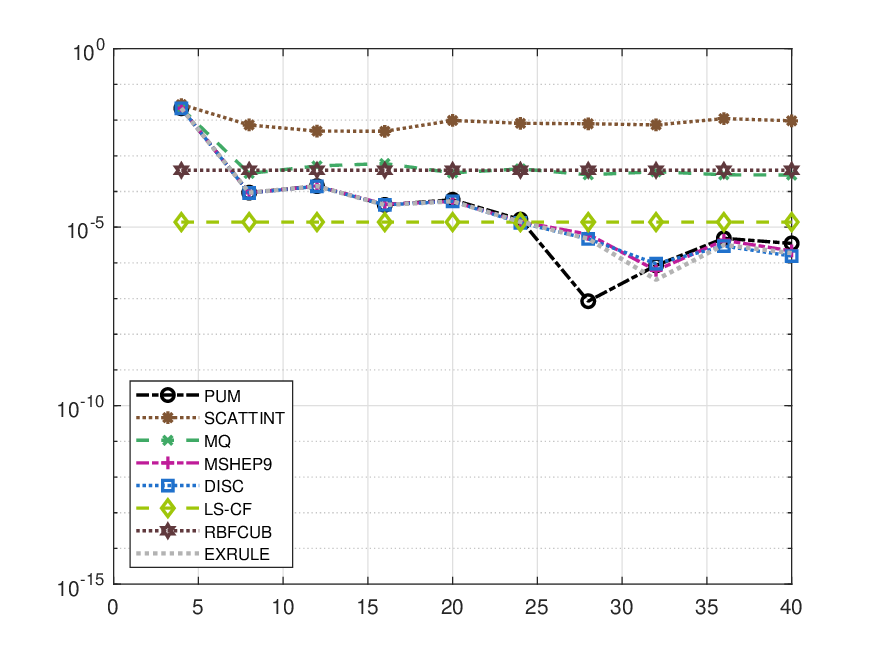}
\includegraphics[scale=0.37]{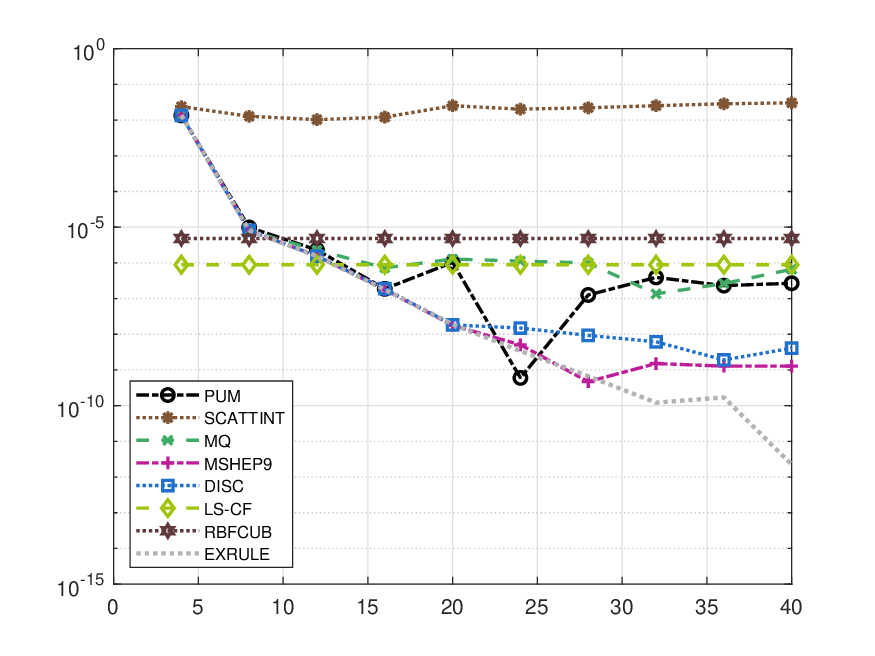}  
\caption{As in Figure \ref{test0} on 800 Halton points.}
\label{test1}
\end{figure}

In Figures \ref{test0}-\ref{test1} we plot the relative cubature errors versus the exactness degree of the underlying (almost-)minimal algebraic rule for the square (such rules have been collected from the literature and inserted in the package \cite{CDRDADTSSV23}). We consider two smooth instances, namely the popular Franke's test function in $\Omega=[0,1]^2$
$$
f_1=\frac{3}{4}\exp(-((9x-2)^2+(9y-2)^2)/4)+\frac{3}{4}\exp(-((9x+1)^2/49+(9y+1)/10))
$$
$$
+\frac{1}{2}\exp(-((9x-7)^2+(9y-3)^2)/4)-\frac{1}{5}\exp(-((9x-4)^2+(9y-7)^2))\;,
$$
and a test function in $\Omega=[-1,1]^2$ from \cite{G21}  
$$
f_2=\frac{1}{(1+x^2)(1+y^2)}\;,
$$
together with two power functions of finite regularity in $\Omega=[0,1]^2$, namely 
$$
f_3=((x-x_0)^2+(y-y_0)^2)^{3/2}\;,\;\;f_4=((x-x_0)^2+(y-y_0)^2)^{7/2}\;,
$$
with $(x_0,y_0)=(0.5,0.5)$, where $f_3\in C^2(\Omega)$ and $f_4\in C^6(\Omega)$. The reference values of the integral have been computed by the same algebraic formula with exactness degree 60.  

\begin{figure}
\centering
\includegraphics[scale=0.3]{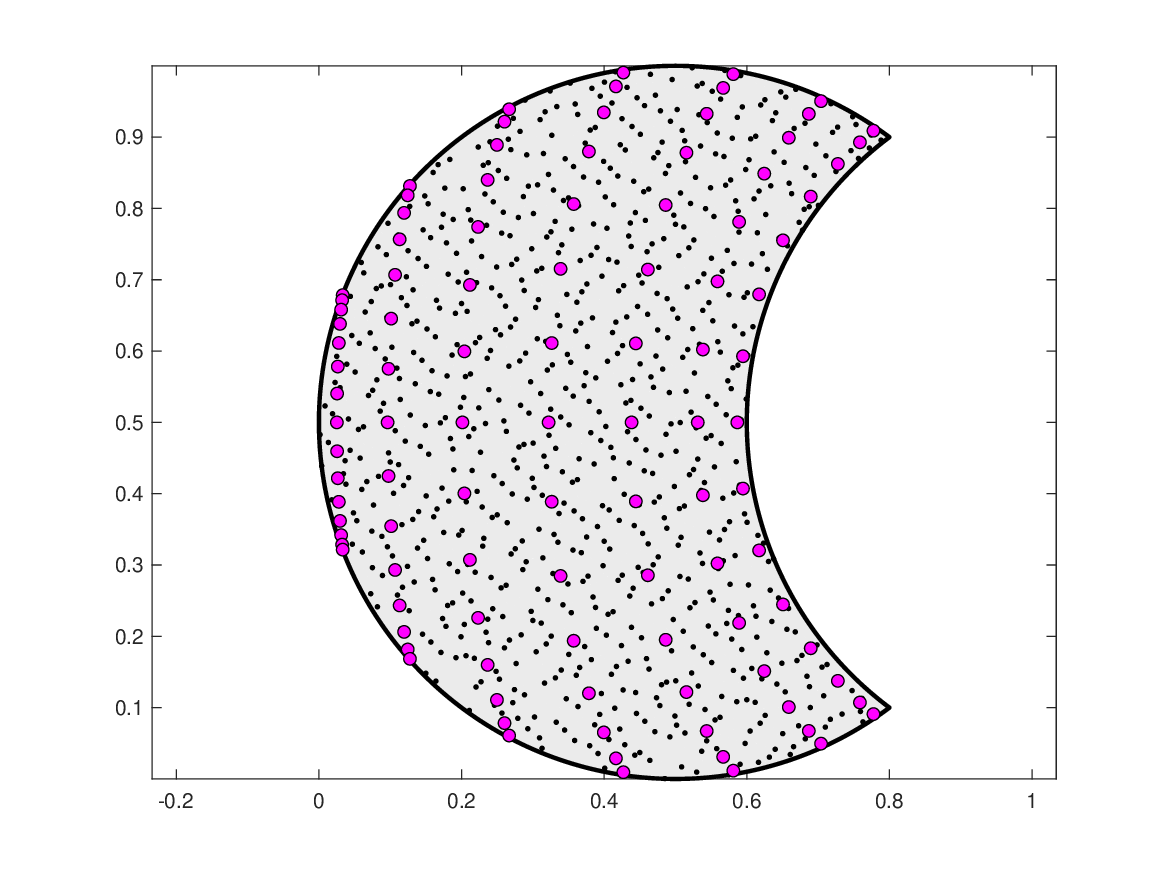}
\includegraphics[scale=0.3]{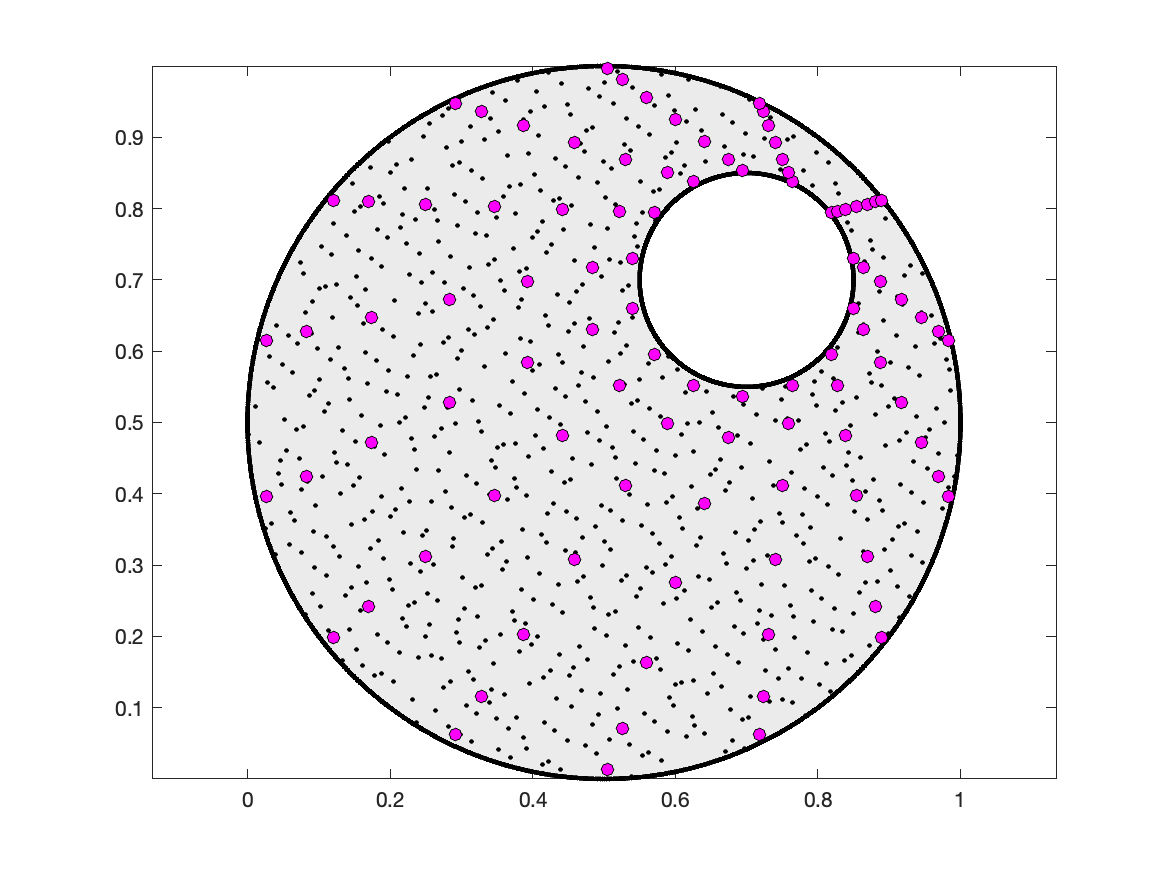}
\caption{Halton points on a circular lune and an asymmetric circular annulus (black dots), and cubature nodes of an algebraic rule of degree 12 (in magenta).}
\label{lune-ann}
\end{figure}

\begin{figure}
\centering
\includegraphics[scale=0.37]{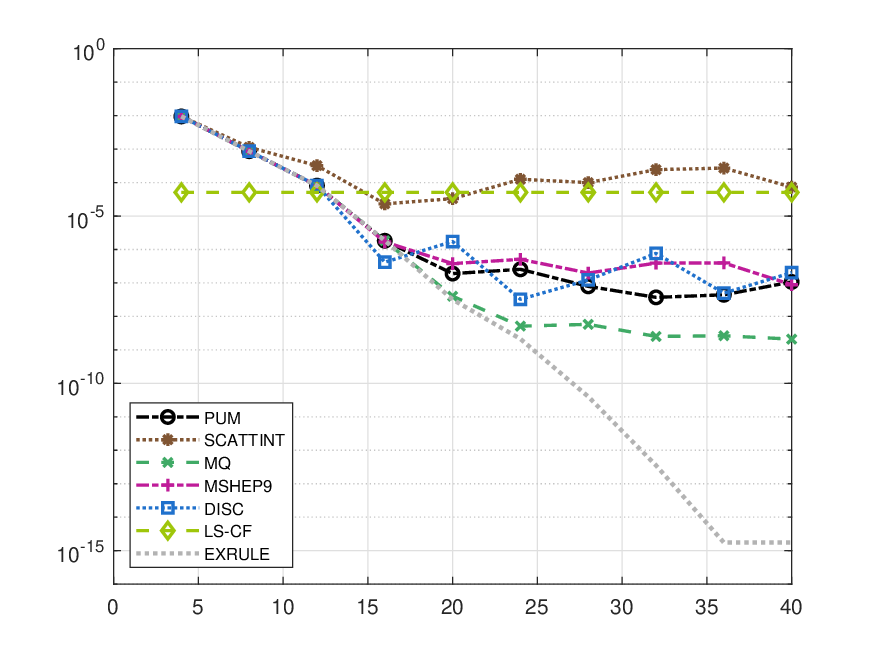}
\includegraphics[scale=0.37]{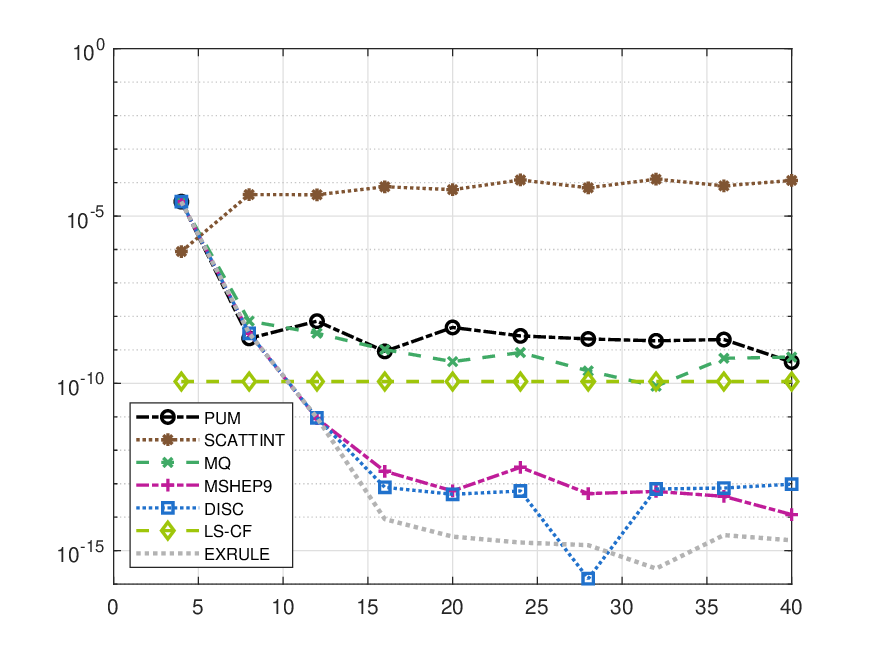}
\includegraphics[scale=0.37]{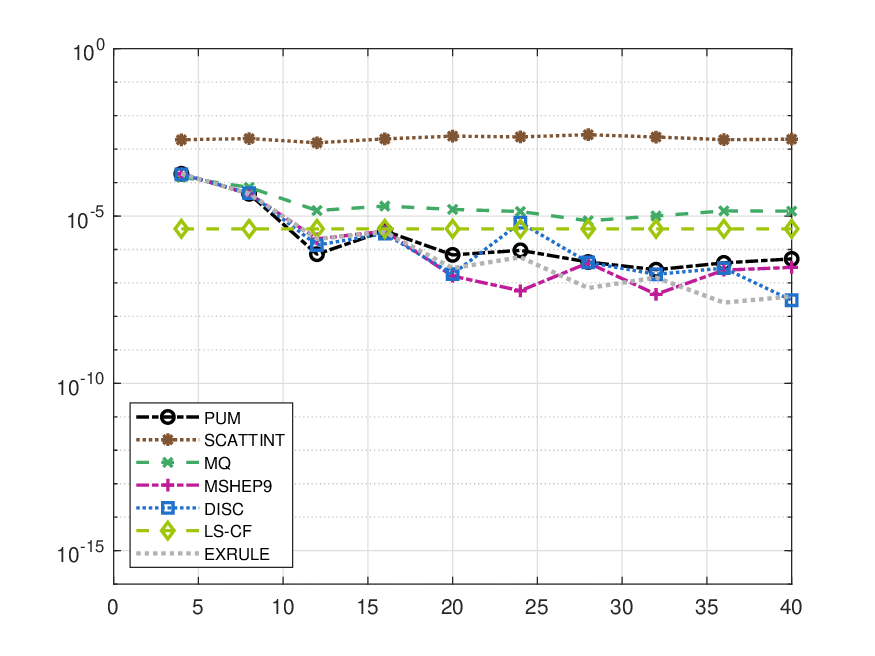}
\includegraphics[scale=0.37]{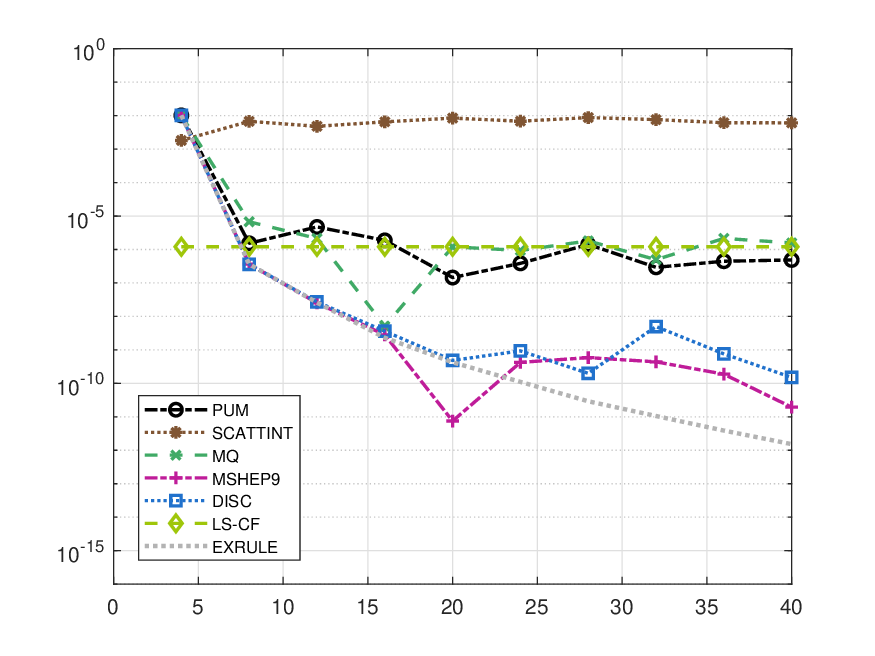}  
\caption{As in Figure \ref{test0} on Halton points of the circular lune of Figure \ref{lune-ann}, extracted 
from 800 Halton points of the bounding box.}
\label{test2}
\end{figure}

\begin{figure}
\centering
\includegraphics[scale=0.37]{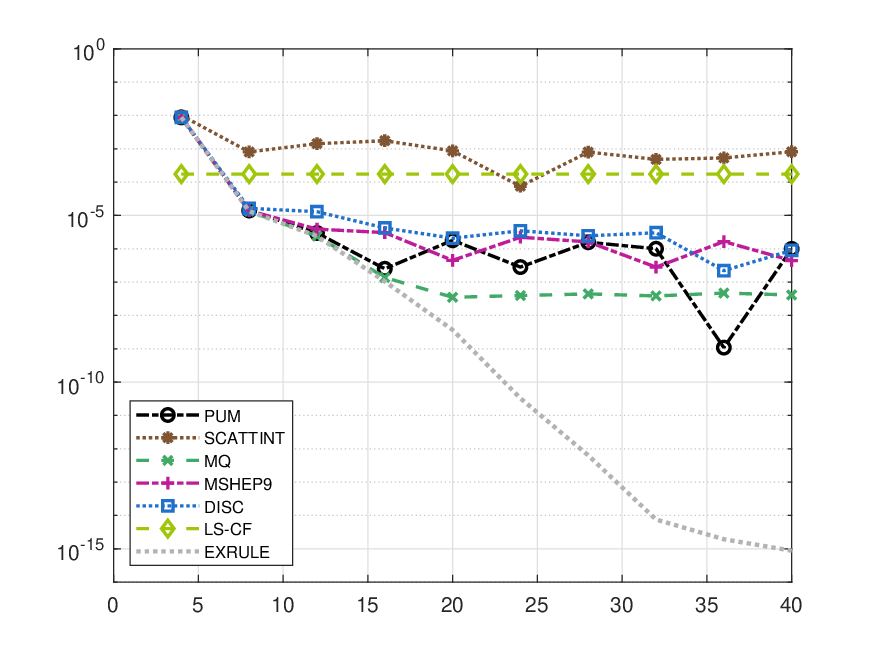}
\includegraphics[scale=0.37]{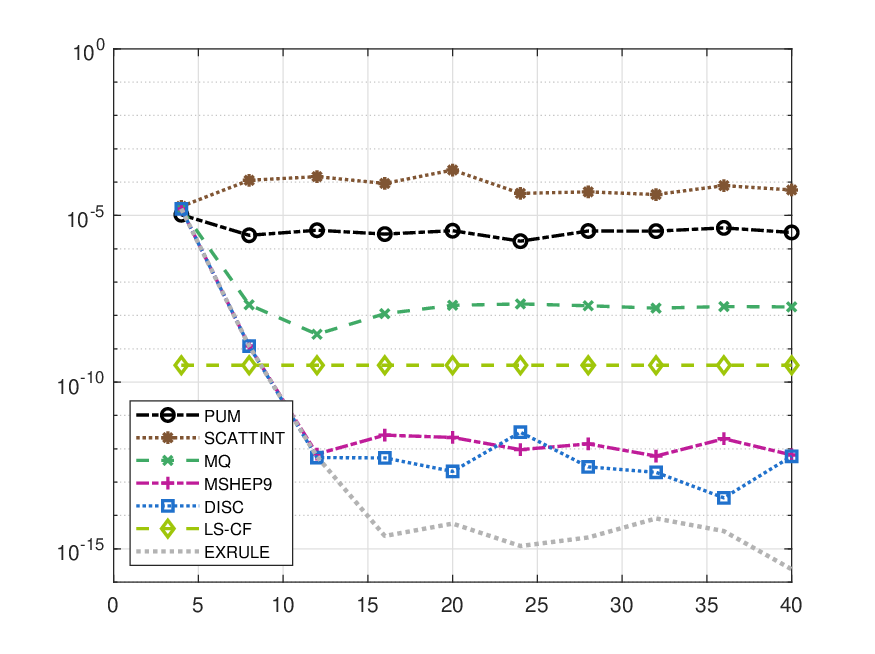}
\includegraphics[scale=0.37]{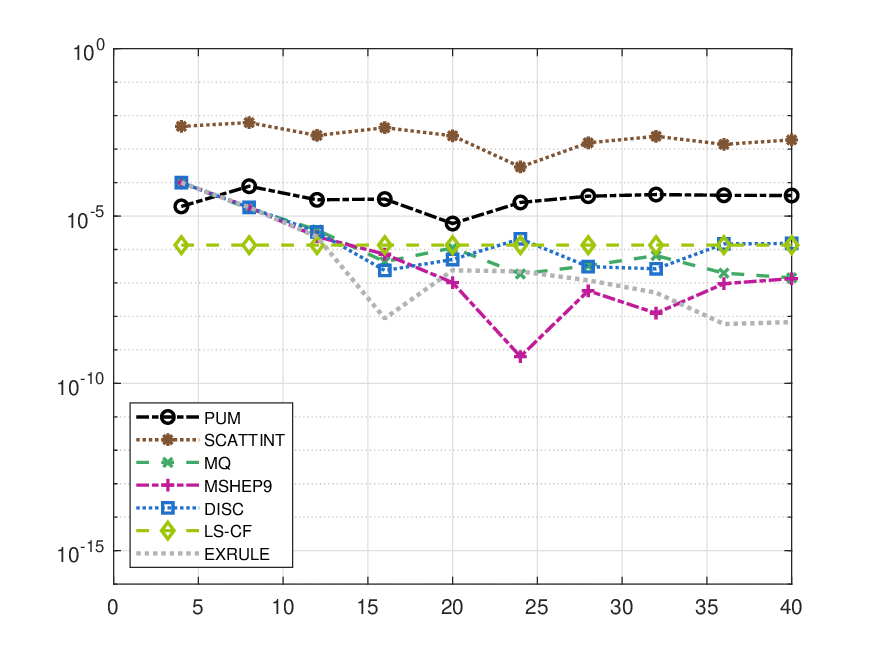}
\includegraphics[scale=0.37]{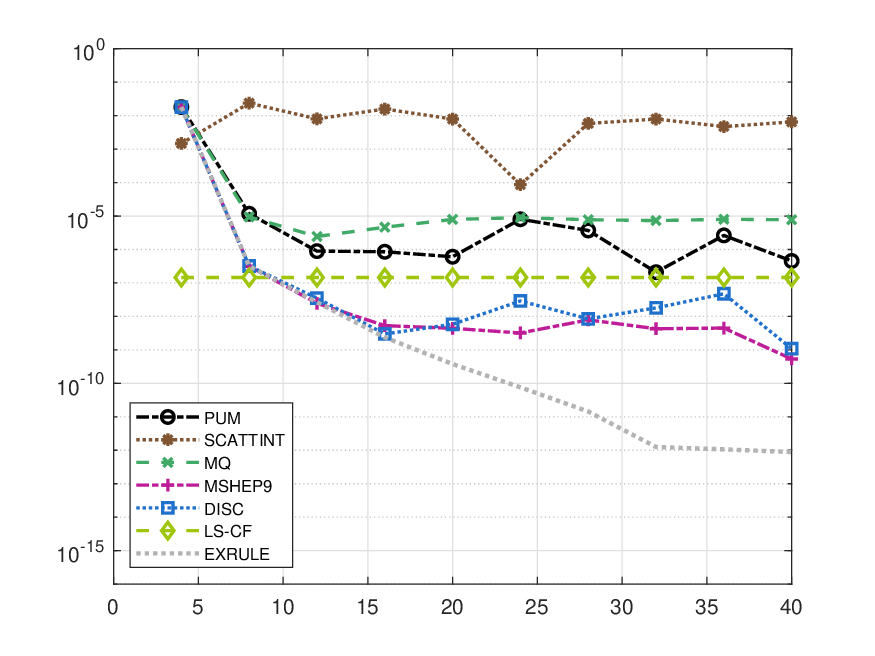}  
\caption{As in Figure \ref{test0} on Halton points of the asymmetric circular annulus of Figure \ref{lune-ann}, extracted 
from 800 Halton points of the bounding box.}
\label{test3}
\end{figure}

The compared methods are Matlab's basic Scattered Interpolation (SCATTINT), adaptive moving interpolation (DISC, \S 2.1), 
adaptive LOOCV RBF interpolation with MultiQuadrics (MQ, \S 2.2), 
BLOOCV RBF-PUM (PUM, \S 2.3), Multinode Shepard method with local interpolation degree 9 (MSHEP9, \S 2.4). The horizontal lines correspond to RBFCUB with MultiQuadrics \cite{CDRSV22} and to LS-CF method \cite{G21}, whereas the grey dotted line to the underlying algebraic rule with exact function values. The choice of MultiQuadrics comes from a number of numerical tests, not reported for brevity, where they have shown to be more accurate than the other RBF considered (Gaussians, IMQ and Wendland $C^2$). 

We see that adaptive moving interpolation and Multinode Shepard follow closely the underlying algebraic rule, until the error stalls around the error generated by  interpolation. It is interesting to observe that higher interpolation errors 
at the boundary (see Figure \ref{Franke}) are there compensated by smaller algebraic cubature weights, preserving the overall cubature accuracy.

Both adaptive moving interpolation and Multinode Shepard at error stalling turn out to be more accurate than Glaubitz's LS-CF method, especially 
with regular integrands where the error goes up to two-three orders of magnitude below. In the comparison with RBFCUB, the situation is similar with the less regular integrands $f_3$ and $f_4$, whereas RBFCUB gives quite accurate results with the smooth integrands $f_1$ and $f_2$. Also PUM method is able to follow on average the underlying algebraic rule until stalling, whereas LOOCV MQ turns out to be quite accurate with the smooth integrands and less accurate, similarly to RBFCUB, with the less regular integrands. 

In order to show the flexibility of the methods, we give also integration examples on two nonstandard domains with curved boundaries, namely a circular lune (non-convex domain) and an asymmetric circular annulus (multiply connected domain),  see Figure \ref{lune-ann}, using the algebraic cubature formulas developed 
in \cite{DFSV13,DFV14} via subperiodic trigonometric Gaussian quadrature. In these cases RBFCUB is not applicable, being restricted to linear polygons. On the other hand, we have modified the code LS-CF, to work with the appropriate polynomial moments on such curved domains (this option is not present in the original package \cite{G21-2}).

Similar considerations with respect to square domains can be done, with adaptive moving interpolation and Multinode Shepard confirming a better accuracy compared 
to LS-CF. It should be stressed, as already observed, that we have here privileged 
maximal accuracy obtainable from not a huge number of fixed data, rather than computational efficiency. Indeed, in the presence of a very high number of scattered points, PUM and LS-CF are faster and could be methods of choice in order to control the computing time.

\vskip0.5cm 
\noindent
{\bf Conclusion and future work.} 
We have constructed cubature methods on scattered data via resampling on the support of known algebraic cubature formulas, by adaptive interpolation of different kinds, such as moving polynomial interpolation at Leja-like points, LOOCV RBF, Multinode Shepard and BLOOCV RBF-PUM. These methods are compared with the most recent approaches to scattered cubature in the literature, namely LOOCV optimized RBF integration over polygons and Glaubitz's LS-CF (based on $\ell^2$ weight minimization under polynomial moment matching conditions).

To our knowledge, this is the first systematic comparison of different approaches for numerical cubature by scattered data, in particular concerning how much accuracy is obtainable from the data. The numerical results confirm that adaptive interpolation applied to known Positive-Interior algebraic cubature formulas can provide valid alternatives to the existing methods, with higher flexibility concerning the shape of integration domains. Extension to three-dimensional integration problems appears natural and will be object of future research.

\vskip0.5cm 
\noindent
{\bf Acknowledgements.} 
Work partially supported by the DOR funds of the University of Padova, by the INdAM-GNCS 2022 Projects ``{\it{Methods and software for multivariate integral models}}'' and ``{\it{Computational methods for kernel-based approximation and its applications}}'', and by the INdAM-GNCS 2023 Project ``{\it{Approximation and multivariate integration, with application to integral equations}}''. 

This research has been accomplished within the RITA ``{\it{Research ITalian network on Approximation}}'', the SIMAI Activity Group ANA\&A, and the UMI Group TAA  ``{\it{Approximation Theory and Applications}}''. 

The first and third authors have been partially supported by the Spoke 1 ``{\it{FutureHPC \& BigData}}'' of ICSC – Centro Nazionale di Ricerca in High-Performance Computing, Big Data and Quantum Computing, funded by European Union - NextGenerationEU.

\end{document}